\documentclass[a4paper, 11pt]{article}

\usepackage[small]{titlesec}
\usepackage{tikz-cd}
\usepackage{booktabs}
\usepackage[T1]{fontenc}
\usepackage[utf8]{inputenc}
\usepackage[english]{babel}
\usepackage[a4paper,top=3.75cm,bottom=3.75cm,left=2.5cm,right=2.5cm]{geometry}
\usepackage{amsfonts}
\usepackage{mathtools}
\usepackage{amsmath}
\usepackage{xcolor}
\mathtoolsset{showonlyrefs}

\usepackage{lmodern}
\usepackage{empheq}
\usepackage{amsthm}
\usepackage{amssymb}
\usepackage{booktabs}
\usepackage{caption}
\usepackage{siunitx}
\usepackage{float}
\usepackage{dsfont}
\usepackage{amssymb}
\usepackage{graphicx}
\usepackage{stmaryrd}

\usepackage{color}
\usepackage{pifont}%
\usepackage{enumerate}
\usepackage{enumitem}
\setlist[description]{leftmargin=\parindent,labelindent=\parindent}
\usepackage{listings} 
\usepackage{amsmath}
\usepackage{amscd}
\usepackage{braket}
\usepackage{tikz}
\usepackage{calc}
\usepackage[hidelinks]{hyperref}
\usepackage{resmes}
\usepackage{mathtools}

\newcommand{\R}{\mathbb{R}}
\newcommand{\sph}{\mathbb{S}}
\newcommand{\K}{\mathcal{K}}

\newcommand{\haus}{\mathcal{H}}
\newcommand{\tra}{\mathrm{tr}}
\newcommand{\Id}{\mathrm{Id}}

\newtheorem{theorem}{Theorem}[section]
\newtheorem*{theorem*}{Theorem}

\newtheorem{proposition}[theorem]{Proposition}
\newtheorem{lemma}[theorem]{Lemma}
\newtheorem{definition}[theorem]{Definition}

\theoremstyle{definition}
\newtheorem{remark}[theorem]{Remark}

\usepackage{pgfplots}
\pgfplotsset{compat=1.17}

\newcommand{\di}{\, \mathrm{d}}

\newcommand{\hbr}{\overline{h}}
\newcommand{\A}{\mathbb{A}}
\newcommand{\B}{\mathbb{B}}
\newcommand{\LB}{\mathcal{L}}
\newcommand{\grad}{\mathrm{grad}}

\usepackage{tikz-cd} 

\newcommand\blfootnote[1]{%
	\begingroup
	\renewcommand\thefootnote{}\footnote{#1}%
	\addtocounter{footnote}{-1}%
	\endgroup
}
\title{A centroaffine perspective on improved Brascamp--Lieb and dimensional Brunn--Minkowski inequalities}
\author{Sotiris Armeniakos, Alexander V. Kolesnikov, and Jacopo Ulivelli }
\date{}

\begin{document}
	\maketitle
	
	\begin{abstract}
		We study dimensional concavity properties of even log-concave measures within the framework of weighted centroaffine differential geometry. We express the relevant local concavity inequalities in terms of a weighted Hilbert--Brunn--Minkowski operator and establish an anisotropic Poincar\'e inequality for odd functions and a strengthened Brascamp--Lieb inequality for even functions. As applications, we obtain quantitative dimensional concavity estimates for a class of geometric measures and measures with pinched Hessian, extending results of Cordero-Erausquin and Rotem for rotationally invariant measures and Eskenazis and Moschidis for the Gaussian measure.
	\end{abstract}

\section{Introduction}

The classical Brunn--Minkowski inequality asserts that, for any two
convex bodies $K,L\subset\mathbb{R}^n$, namely compact convex sets with
nonempty interior, and every $t\in[0,1]$,
\begin{equation}\label{eq:BM}
    |(1-t)K+tL|^{\frac1n}
    \geq (1-t)|K|^{\frac1n}+t|L|^{\frac1n}.
\end{equation}
Here, $|\cdot|$ denotes the Lebesgue measure. Moreover, we denote by
$\mathcal K^n$ the family of convex bodies in $\mathbb R^n$.
Considerable attention has been devoted to extending and strengthening
\eqref{eq:BM}, particularly under symmetry assumptions. One direction
consists of replacing Minkowski addition by the $L_p$-addition of convex
bodies. Of particular interest is the range $0\leq p<1$, where the
$L_p$-Brunn--Minkowski conjectures assert that
\begin{equation}\label{eq:L_p-BM}
    |(1-t)K+_p tL|
    \geq
    \left(
        (1-t)|K|^{\frac pn}+t|L|^{\frac pn}
    \right)^{\frac np},
\end{equation}
for all origin-symmetric $K,L\in\mathcal K^n$ and every $t\in[0,1]$.
We defer the precise definition of $L_p$-addition to Section~2 and
note that, in the chosen range $0 \leq p <1$, \eqref{eq:L_p-BM} strengthens \eqref{eq:BM}.
The limiting case $p=0$ is the celebrated log-Brunn--Minkowski
conjecture of B\"or\"oczky, Lutwak, Yang and
Zhang~\cite{boroczkyLYZ2013}, which asserts that
\begin{equation}\label{eq:log-BM}
    |(1-t)K+_0tL|\geq |K|^{1-t}|L|^t.\notag
\end{equation}

\medskip

An approach to these questions is provided by the spectral properties
of the Hilbert--Brunn--Minkowski operator, which we briefly describe.
Let $\mathcal K^{2,\alpha}_{e,+}$ denote the class of origin-symmetric
convex bodies with $C^{2,\alpha}$ boundary and strictly positive
principal curvatures. For $K\in\mathcal K^{2,\alpha}_{e,+}$, let $h_K$
be its support function and define the tensor
\[
    Q_K=\nabla^2_{\mathbb S^{n-1}}h_K+h_K\operatorname{Id}
\]
on the tangent bundle of $\mathbb S^{n-1}$. Here, $\nabla_{\sph^{n-1}}$ denotes the spherical gradient. 
With the normalization used here, the Hilbert--Brunn--Minkowski
operator acts on $z\in C^2(\mathbb S^{n-1})$ by
\begin{equation}\label{eq:Hilbert_operator}
    \Delta_Kz
    =(Q_K^{-1})^{ij}
    \left(
        h_K\nabla^2_{ij}z+(h_K)_iz_j+(h_K)_jz_i
    \right),
\end{equation}
where the standard covariant derivative on $\sph^{n-1}$ is considered and repeated indices are summed.
Kolesnikov and E.~Milman~\cite{Kol_Mil_p, Mil_centroaffine} related local
$L_p$-Brunn--Minkowski inequalities to lower bounds on the first
nonzero even eigenvalue of $-\Delta_K$. In particular, they
established the local inequalities for every smooth origin-symmetric
convex body in a dimension-dependent range of parameters $p<1$
sufficiently close to $1$. E.~Milman~\cite{Mil_centroaffine} subsequently identified $\Delta_K$ with the Laplace--Beltrami operator on the boundary of $K$ (which we denote by $\partial K$), associated with a suitable Riemannian metric and a connection distinct from the Levi-Civita connection of that metric. Doing so, the treatment in \cite{Mil_centroaffine} placed the problem in the context of so-called centroaffine differential geometry and led to results under curvature pinching conditions. The relation between this geometry and the Hessian metric of the Minkowski functional was further developed by Ivaki and E.~Milman~\cite{Ivaki_Milman_p}. Similar results under pinching were later obtained in connection to improved Brascamp--Lieb inequalities by Kolesnikov, Livshyts and Rotem \cite{Kol_liv_Rot}.

\medskip

Another direction concerns replacing the Lebesgue measure by more general measures. For a measure $\mu$ with density $e^{-u}$, where $u$ is convex, the Pr\'ekopa--Leindler inequality asserts that \begin{equation}\label{eq:log-concavity_measure}
    \mu((1-t)K+tL)\geq\mu(K)^{1-t}\mu(L)^t
\end{equation} 
for all convex bodies $K,L$ and every $t\in[0,1]$. Equivalently (compare the work of Borell \cite{Borell}), a Borel measure $\mu$ satisfying \eqref{eq:log-concavity_measure} is said to be log-concave. For the Lebesgue measure, it is immediate via homogeneity that this multiplicative adimensional inequality is equivalent to \eqref{eq:BM}. The same procedure evidently breaks for general log-concave measures. Improving the concavity of log-concave measures is a research topic that has produced a significant amount of work in the study of log-concave measures, the main example being the dimensional Brunn--Minkowski conjecture: It postulates that the volume in \eqref{eq:BM} can be replaced with an arbitrary even log-concave measure when the convex sets are assumed to be origin symmetric. As a natural question, the weighted and $L_p$ extensions meet in the study of inequalities
of the form
\begin{equation}\label{eq:weighted_Lp_BM_intro}
    \mu((1-t)K+_p tL)
    \geq
    \left(
        (1-t)\mu(K)^{\frac pn}
        +t\mu(L)^{\frac pn}
    \right)^{\frac np},
\end{equation}
for origin-symmetric $K, L \in \K^n$, and even log-concave measures $\mu$, with $p=0$ interpreted in the limiting sense. The choice $p=1$ gives the dimensional Brunn--Minkowski conjecture, while $p=0$ gives the weighted log-Brunn--Minkowski conjecture. These inequalities were studied within a common $(p,q)$-concavity framework by Hosle,
Kolesnikov and Livshyts \cite{hos_kol_liv}.

\medskip

For the choice of $\mu$ as the Gaussian measure and $p=1$, Gardner and Zvavitch~\cite{gardner2010gaussian} originally proposed the dimensional inequality \eqref{eq:weighted_Lp_BM_intro} under the weaker assumption that the convex bodies contain the origin. Nayar and Tkocz~\cite{nayar2013note} disproved this formulation, leaving the origin-symmetric case as the principal question. The corresponding dimensional problem for general even log-concave measures was investigated by Colesanti, Livshyts and Marsiglietti~\cite{ColLivMar_stability}, who established local results near the Euclidean ball for rotationally invariant log-concave measures.

We briefly comment on the usual approach for this weighted setting, using a Bochner formula due to Kolesnikov and E.~Milman~\cite{Kol_Mil_BL,Kol_Mil_Po} and Kolesnikov and Livshyts \cite{Kol_Liv_GZ}, for these concavity inequalities it is sufficient to establish coercivity estimates for quantities of the form

\[\int_{K} \|D^2f\|_{HS}+\langle (D^2u)Df,Df \rangle \di \mu \]
for every origin symmetric convex set $K$ and even function $f$. Using this approach, Eskenazis and Moschidis~\cite{Esk_Mos_Gauss} subsequently proved the Gaussian inequality with the conjectured exponent $1/n$ for all origin-symmetric convex bodies. Further extensions have since been obtained. Han and Zhu~\cite{HanZhuPoincare} replaced origin symmetry by the assumption that the Gaussian measures conditioned on the two convex sets have the same barycenter, and obtained an inequality with explicit barycenter terms in the general case. Recently Milman, Nakamura and Tsuji~\cite{mil_nak_tsu} established a heterovariant Gaussian Brunn--Minkowski inequality, involving Gaussian measures with different covariance matrices subject to a matrix compatibility condition. Their result contains the dimensional Gaussian inequality and connects it with strengthened Gaussian dilation inequalities, which relates to the (B) conjecture, which we now state. The Gaussian statement was proposed by Banaszczyk and popularized by Lata{\l}a~\cite{Latala2002ICM}; its extension to even log-concave measures asks whether, for every origin-symmetric convex body $K$, the function
\begin{equation}\label{eq:B_conjecture}
    t\longmapsto\mu(e^tK)
\end{equation}
is log-concave on $\mathbb R$. Cordero-Erausquin, Fradelizi and Maurey~\cite{Cor_Fra_Mau} proved the Gaussian case. We also remark that the (B) property holds for a certain class of Gaussian mixtures \cite{EskeNayTkoczMixtures}. Saroglou~\cite{Sar} showed that the log-Brunn--Minkowski conjecture for the Lebesgue measure implies its weighted version for even log-concave measures and, consequently, the (B) conjecture. Cordero-Erausquin and Rotem~\cite{rot_B} established both the (B) property and the dimensional Brunn--Minkowski inequality for rotationally invariant log-concave measures, as well as for a broader class of measures with radial densities. 

\medskip

The general scheme is to prove an improved Poincar\'e inequality for odd functions, then apply this estimate to the coordinate functions of a gradient of an even function and conclude some coercivity estimate. In the general case that the measure is assumed to be merely log-concave this was combined with delicate gradient estimates to produce a weaker form of \eqref{eq:weighted_Lp_BM_intro}, in an approach put forth by Livshyts \cite{Liv_Dimensional_BM} and extended by Eskenazis, Giannopoulos and Tziotziou~\cite{eskenazis2026functionalperimeterdimensionalbrunnminkowski} and Yang \cite{yang2026uniformbounddimensionalbrunnminkowski}. 

\medskip
\paragraph{Main results.}
The purpose of this work is to connect the spectral and centroaffine approach to the log-Brunn--Minkowski problem with the functional inequalities developed for dimensional concavity and the (B) conjecture. A similar path has been recently explored also by Hu and Ivaki \cite{Ivaki_Hu_weighted_centroaffine}.

\medskip

We first express the weighted local concavity inequalities in terms of a spectral minimization problem for an appropriate operator defined on the boundary of convex sets. More precisely, given $K \in \mathcal K^{2,\alpha}_{e,+}$ and a log-concave measure $\mu$ with density of class $C^1$, we define the weighted Hilbert--Brunn--Minkowski operator as
\begin{equation}\label{eq:wHBM}
    \mathcal{D}_{K,\mu}f \coloneqq -\Delta_{K,\mu}f + \langle \nabla u,x \rangle f. 
\end{equation}
Here, $\Delta_{K,\mu}$ denotes the weighted centroaffine Laplace--Beltrami operator, which will be defined properly in Section~2. Moreover, $\langle \cdot, \cdot \rangle$ denotes the standard Euclidean scalar product in $\R^n$. We use the same notation for its restriction to embedded submanifolds of $\R^n$. The weighted cone-volume measure $V_{K,\mu}$ on $\partial K$ with respect to $\mu$ is defined as the measure on $\partial K$ with density
\[  \di V_{K,\mu}(x)=\frac{\langle x, \nu_{K}(x)\rangle}{n\mu(K)}\di \mathcal{H}^{n-1}.\]
We note that pushing \eqref{eq:wHBM} on the sphere through the Gauss map in the unweighted setting yields precisely the usual Hilbert--Brunn-Minkowski operator \eqref{eq:Hilbert_operator}. 

\medskip

With these definitions at hand, the relevant functional for us is the following: 
\[
    \mathcal J_{K,\mu}(f)
    \coloneqq
    \frac{
        \displaystyle\int_{\partial K}
        (\mathcal D_{K,\mu}f)f\,\mathrm dV_{K,\mu}
    }{
        \int_{\partial K} f^2 \di V_{K,\mu} - \left(\int_{\partial K} f \di V_{K,\mu}\right)^2
    }.
\]
We remark that whenever $\di \mu(x) = e^{-u(x)}\di x$ and $u$ is strictly convex and even, the denominator is never zero unless $f \equiv 0$. Our first result records the common spectral formulation of the preceding concavity properties.
\begin{theorem}\label{thm:common_spectral_framework}
    Let $\mu$ be an even log-concave measure on $\mathbb R^n$
    with density $e^{-u}$, where $u\in C^1(\mathbb R^n)$ is
    strictly convex. Then the following characterizations hold:
    \begin{enumerate}
        \item For $p\in[0,1]$, inequality
        \eqref{eq:weighted_Lp_BM_intro} holds for all
        origin-symmetric convex bodies $K,L$ and every
        $t\in[0,1]$ if and only if
        \begin{equation}\label{eq:weighted_eigenvalue}
            \inf_{\substack{
                f\in C_e^2(\partial K)\\
                \mathrm{Var}_{V_{K,\mu}}(f)>0
            }}
            \mathcal J_{K,\mu}(f)\geq n-p
        \end{equation}
        for every $K\in\mathcal K^{2,\alpha}_{e,+}$.
        The case $p=0$ is understood in the limiting sense.

        \item The measure $\mu$ satisfies the $(B)$ inequality
        for all origin-symmetric convex bodies if and only if
        \[
            \inf_{\substack{f\equiv c\neq0}}
            \mathcal J_{K,\mu}(f)\geq n
        \]
        for every $K\in\mathcal K^{2,\alpha}_{e,+}$.
    \end{enumerate}
\end{theorem}\noindent
We note that in the case of the Lebesgue measure, the weighted cone-volume measure $V_{K,\mu}$ reduces to the usual normalized cone-volume measure and \eqref{eq:weighted_eigenvalue} reduces to the first even eigenvalue of the Hilbert--Brunn--Minkowski operator.

\medskip

Fix now a smooth origin-symmetric convex body $K$, and set
\[
    \Phi(x)=\frac{\|x\|_K^2}{2},
\]
where $\|\cdot\|_K$ denotes its Minkowski functional. Following the ideas of \cite{rot_B}, we were able to extend their central Poincar\'e inequality for odd functions to the anisotropic setting. Let us remark that, during the preparation of this work, a preprint by Cordero-Erausquin, Livshyts and Rotem \cite{cor-er_liv_rot} appeared where the following result was obtained independently. The proof strategy is the same up to minor differences, one of them being that we use a weighted centroaffine formalism. This shall be of no surprise, as the blueprint is the one laid down in the substantial work by Cordero-Erausquin and Rotem \cite{rot_B}. Nonetheless, in this manuscript we manage to further leverage the ideas from \cite{rot_B} combining them with the approach of \cite{Esk_Mos_Gauss}. In particular, we obtain the conclusions expected in \cite[p. 9]{cor-er_liv_rot} as a consequence of the following result. We present these implications in the remainder of this manuscript.
\begin{theorem}\label{thm:intro_anisotropic_odd}
    Let $K$ be a smooth origin-symmetric convex body with strictly
    positive principal curvatures, let
    $w\in C^1((0,\infty))$ be increasing, and let
    $V:\mathbb R^n\to\mathbb R$ be even and convex.
    Consider the finite measure
    \[
        \mathrm d\eta(x)
        =e^{-V(x)-w(\|x\|_K)}\,\mathrm dx.
    \]
    Then, for every odd $\varphi\in C^1(\mathbb R^n)$ for which
    the integrals are finite,
    \begin{equation}\label{eq:intro_odd_poincare}
        \int_{\mathbb R^n}
            \frac{w'(\|x\|_K)}{\|x\|_K}
            \varphi^2\,\mathrm d\eta
        \leq
        \int_{\mathbb R^n}
            \left\langle
                (D^2\Phi)^{-1}\nabla\varphi,
                \nabla\varphi
            \right\rangle\mathrm d\eta.
    \end{equation}
\end{theorem}

To pass from the above to an improved Brascamp--Lieb inequality for even functions, we impose the pinching condition
\begin{equation}\label{eq:intro_gauge_pinching}
    \alpha\operatorname{Id}
    \leq D^2\Phi(x)
    \leq\beta\operatorname{Id},
    \qquad x\neq0,
    \qquad 0<\alpha\leq\beta<\infty,
\end{equation}
and write $\gamma=\alpha/\beta$. Adapting the approach from \cite{rot_B}, we obtain the following consequence of Theorem \ref{thm:intro_anisotropic_odd}.
\begin{theorem}\label{thm:intro_anisotropic_even}
    Let $K$ be a smooth origin-symmetric convex body with strictly
    positive principal curvatures satisfying
    \eqref{eq:intro_gauge_pinching}, and let
    $w\in C^2((0,\infty))$ be increasing and satisfy
    \begin{equation}\label{eq:intro_profile_condition}
        w''(s)+\gamma\frac{w'(s)}s\geq0,
        \qquad s>0.
    \end{equation}
    Let $V\in C^2(\mathbb R^n)$ be even and convex, and suppose that
    \[
        \mathrm d\eta(x)=e^{-V(x)-W(x)}\,\mathrm dx,
        \qquad W(x)=w(\|x\|_K),
    \]
    is a probability measure. Then, for every even
    $\varphi\in C^1(\mathbb R^n)$,
    \begin{equation}\label{eq:intro_even_brascamp_lieb}
        \operatorname{Var}_{\eta}(\varphi)
        \leq
        \int_{\mathbb R^n}
        \left\langle
            \left(
                D^2W+
                \gamma\frac{w'(\|x\|_K)}{\|x\|_K}
                D^2\Phi
            \right)^{-1}
            \nabla\varphi,\nabla\varphi
        \right\rangle\mathrm d\eta.
    \end{equation}
    When the matrix is singular, the inverse quadratic form is
    understood as the associated dual quadratic form.
\end{theorem}\noindent 
Here $\operatorname{Var}_\eta(\varphi)
=\int\varphi^2\,\mathrm d\eta-(\int\varphi\,\mathrm d\eta)^2$.
When $K$ is the Euclidean unit ball,
$D^2\Phi$ is just the identity matrix $\operatorname{Id}$ and one can take $\gamma=1$, recovering the corresponding inequalities of Cordero-Erausquin and Rotem~\cite{rot_B}. 

\medskip

Using the $L_2$ method and the above inequalities we can conclude a dimensional Brunn--Minkowski inequality for this class of measures.
\begin{theorem}\label{thm:intro_dimensional_gauge}
    Let $K$ be a smooth origin-symmetric convex body with strictly
    positive principal curvatures, and suppose that
    $\Phi=\|\cdot\|_K^2/2$ satisfies
    \eqref{eq:intro_gauge_pinching}. Let $\gamma=\alpha/\beta$,
    and let $w\in C^2((0,\infty))$ be increasing and satisfy
    \eqref{eq:intro_profile_condition}. Suppose that
    \[
        \mathrm d\mu(x)
        =e^{-w(\|x\|_K)}\,\mathrm dx
    \]
    defines a finite measure. Then, for every pair of
    origin-symmetric convex bodies $A,B\subset\mathbb R^n$
    and every $t\in[0,1]$,
    \begin{equation}\label{eq:intro_dimensional_gauge}
        \mu((1-t)A+tB)^{\frac\gamma n}
        \geq
        (1-t)\mu(A)^{\frac\gamma n}
        +t\mu(B)^{\frac\gamma n}.
    \end{equation}
\end{theorem} \noindent
We point out that inequality \eqref{eq:intro_dimensional_gauge} constitutes a quantitative version of the result of \cite{rot_B}, since the parameter $\gamma$ in the pinching above quantifies how far the set $K$ is from a ball, in centroaffine terms. We remark that some of the measures above need not necessarily be log-concave. 

Finally, we show that a similar argument, independent of the centroaffine machinery we develop, works for general measures with a suitable pinching on the potential, leading to the following result.
\begin{theorem}\label{thm:intro_dimensional_pinched}
    Let $V\in C^2(\mathbb R^n)$ be even and convex, and suppose that
    \[
        k_1\operatorname{Id}\leq D^2V,
        \qquad
        \Delta V\leq k_2n
    \]
    for some $0<k_1\leq k_2<\infty$.
    Set $r=k_1/k_2$ and let
    $\mathrm d\mu=e^{-V}\,\mathrm dx$.
    Then, for every pair of origin-symmetric convex bodies
    $A,B\subset\mathbb R^n$ and every $t\in[0,1]$,
    \begin{equation}\label{eq:intro_dimensional_pinched}
        \mu((1-t)A+tB)^{\frac rn}
        \geq
        (1-t)\mu(A)^{\frac rn}
        +t\mu(B)^{\frac rn}.
    \end{equation}
\end{theorem}\noindent
Theorem \ref{thm:intro_dimensional_pinched} can be understood as an extension of the result of Eskenazis and Moschidis \cite{Esk_Mos_Gauss}, since for the Gaussian $k_{1}=k_{2}=1$. In particular, notice that Theorem \ref{thm:intro_dimensional_pinched} improves the previous dimensional inequalities from \cite{Kol_Liv_GZ} and \cite{Armen_Uli} under pinching conditions.

\medskip

\paragraph{Structure of the manuscript.} The paper is organized as follows. Section~2 collects the geometric notation and the properties of the weighted centroaffine Laplacian. Section~3 develops the weighted Hilbert--Brunn--Minkowski framework and discusses its relation to local concavity inequalities and spectral comparison questions. Section~4 establishes the anisotropic Poincar\'e and strengthened Brascamp--Lieb inequalities. Section~5 derives the dimensional concavity estimates for measures associated with a Minkowski functional and for general pinched potentials.

\paragraph{Further context for the results.} As previously mentioned in this introduction, both the dimensional Brunn--Minkowski conjecture and the (B) conjecture are implied by the logarithmic Brunn--Minkowski conjecture. An AI-generated solution to this major problem was released by OpenAI at the beginning of October 2026. The elaboration of this manuscript and its content dates back to prior to the aforementioned release, and our proofs follow the path laid out by the various works mentioned in the introduction. We believe our results to be of independent interest for a deeper understanding of the aforementioned conjectures. Moreover, several of our statements involve measures beyond log-concave ones and, therefore, are not implied by the logarithmic Brunn--Minkowski conjecture.

\section{Preliminaries}

Throughout this work, let $K\subset\R^n$ be a compact convex set
whose boundary is a $C^2$ manifold with bounded and strictly positive principal curvatures. We call a set of this type a smooth convex body and assume that the origin lies in its interior. We denote by $\nu_K:\partial K\to\sph^{n-1}$ the corresponding Gauss map and by ${\rm I\! I}$ its second fundamental form. Identifying the latter with the associated shape operator, we have $\di\nu_K={\rm I\! I}$. We refer to \cite[Section~2.5]{Schneider_CB} for further details on the differential geometry of smooth convex bodies.

Let $h_K:\sph^{n-1}\to\R, h_K(\xi)={\sup}_{x \in K}\langle x, \xi \rangle$, denote the support function of $K$, and
write $\hbr_K=h_K\circ\nu_K$, so that $\hbr_K$ is defined on $\partial K$. Following
\cite{Kol_Mil_p,Mil_centroaffine}, we equip $\partial K$ with the
centroaffine metric
\[
    g_K=\frac{{\rm I\! I}}{\hbr_K}.
\]
For computational purposes, we work directly on $\partial K$.
The corresponding spherical formulation is obtained through the
Gauss map; see also
\cite{Ivaki_Milman_p,Ivaki_Hu_weighted_centroaffine}. For convex bodies $K,L\in\mathcal K^n$ containing the origin in their interiors, $p\in(0,1)$ and $t\in[0,1]$, we define their $L_p$-sum by
\begin{equation}\label{eq:Lp_sum}
    (1-t)K+_p tL
    \coloneqq
    \bigcap_{\theta\in\sph^{n-1}}
    \left\{
        x\in\R^n:
        \langle x,\theta\rangle
        \leq
        \left(
            (1-t)h_K(\theta)^p+t h_L(\theta)^p
        \right)^{1/p}
    \right\}.
\end{equation}
This is the Wulff shape associated with the function $((1-t)h_K^p+t h_L^p)^{1/p}$, which need not itself be a support function. The limiting case $p=0$ is defined by replacing this function with $h_K^{1-t}h_L^t$, while $p=1$ recovers the usual Minkowski combination $(1-t)K+tL$. 

Let $u\in C^1(\R^n)$ be convex, and let $\mu$ be the log-concave measure defined by
\[
    \di\mu(x)=e^{-u(x)}\di x.
\]
With a slight abuse of notation, we also use $\mu$ to denote the surface measure $e^{-u}\di\haus^{n-1}$ on $\partial K$. We define the $\mu$-cone-volume measure of $K$ by
\[
    \di V_{K,\mu}(x)
    =\frac{\hbr_K(x)e^{-u(x)}}{n\mu(K)}
    \di\haus^{n-1}(x).
\]
For the Lebesgue measure, this agrees with the normalized cone-volume measure. In general, however, $V_{K,\mu}$ need not be a probability measure, a distinction that will be relevant in what follows.

Building on \cite{Kol_Mil_p,Mil_centroaffine,Ivaki_Milman_p}, we define, for $f\in C^2(\partial K)$,
\[
\begin{aligned}
    \Delta_{K,\mu}f
    ={}
    \hbr_K\nabla_{\partial K}\cdot
        \left({\rm I\! I}^{-1}\nabla_{\partial K}f\right)
    +2\left\langle
        {\rm I\! I}^{-1}\nabla_{\partial K}f,
        \nabla_{\partial K}\hbr_K
    \right\rangle-\hbr_K\left\langle
        {\rm I\! I}^{-1}\nabla_{\partial K}f,
        \nabla_{\partial K}u
    \right\rangle.
\end{aligned}
\]
When $u\equiv0$, this reduces to the Hilbert--Brunn--Minkowski operator with the normalization of \cite{Mil_centroaffine}. Equivalently, $\Delta_{K,\mu}$ is the weighted Laplace--Beltrami operator associated with $(\partial K,g_K,V_{K,\mu})$. Indeed, the corresponding metric gradient is
\[
    {\rm grad}_{g_K}f
    =g_K^{-1}\nabla_{\partial K}f
    =\hbr_K{\rm I\! I}^{-1}\nabla_{\partial K}f,
\]
and integration by parts gives, for every $f,g\in C^2(\partial K)$, \begin{align}\label{eq:IBP_w_aff}
    \int_{\partial K}
        (-\Delta_{K,\mu}f)g\,\di V_{K,\mu}
    &=
    \int_{\partial K}
        g_K\left(
            {\rm grad}_{g_K}f,
            {\rm grad}_{g_K}g
        \right)\di V_{K,\mu}\\
    &=
    \int_{\partial K}
        \hbr_K\left\langle
            {\rm I\! I}^{-1}\nabla_{\partial K}f,
            \nabla_{\partial K}g
        \right\rangle\di V_{K,\mu}.
        \notag
\end{align}

\section{On the weighted Hilbert operator}

In this section, we define the weighted analogue of the Hilbert
operator, following \cite{Kol_Mil_p}, and study its connection
to local formulations of global concavity inequalities and their
spectral counterparts. Our starting point is the $(p,q)$-concavity
inequality
\begin{equation}\label{eq:global (p,q)-Version}
    \mu((1-t)K+_{p}tL)^{\frac{q}{n}}
    \geq (1-t)\mu(K)^{\frac{q}{n}}
    +t\mu(L)^{\frac{q}{n}},
\end{equation}
where $0\leq q\leq p\leq1$, $K,L$ are origin-symmetric convex
bodies, and $\mu$ is an even log-concave probability measure
on $\R^n$ with density
\[
    \di\mu=e^{-u}\di x,
\]
for an even, convex, lower semicontinuous function $u:\R^n\to\R\cup\{+\infty\}$. The case $q=0$ is understood in the logarithmic sense. For the differential and spectral statements below, we additionally assume that $u$ is finite and belongs to $C^2(\R^n)$. 

As established in \cite[Lemma~2.1]{hos_kol_liv}, the validity of \eqref{eq:global (p,q)-Version} for all origin-symmetric convex bodies implies the following infinitesimal inequality: for every origin-symmetric $K\in\mathcal K^n_{2,+}$ and every even $f\in C^1(\partial K)$,
\begin{equation}\label{eq:local (p,q)-Version}
\begin{aligned}
    &\int_{\partial K}
    \left(
        \left\langle
            {\rm I\! I}^{-1}\nabla_{\partial K}f,
            \nabla_{\partial K}f
        \right\rangle
        -{\rm H}_\mu f^2
        -(1-p)\frac{f^2}{\hbr_K}
    \right)\di\mu+\frac{n-q}{n\mu(K)}
        \left(\int_{\partial K}f\di\mu\right)^2
    \geq0,
\end{aligned}
\end{equation}
where
${\rm H}_\mu=\mathrm{tr}({\rm I\! I})-\langle\nabla u,\nu_K\rangle$. As in the preliminaries, $\di\mu$ on $\partial K$ denotes $e^{-u}\di\haus^{n-1}$. 

The converse implication is more delicate. In the unweighted setting, local-to-global arguments were developed using strongly isomorphic polytopes in \cite{Putter2021} and a PDE approach in \cite{LocalGlobalPDE}. In the weighted setting, \cite[Theorem~3.1]{hos_kol_liv} shows that, for $0\leq q\leq p<1$, the validity of \eqref{eq:local (p,q)-Version} for all such $K$ and $f$ implies \eqref{eq:global (p,q)-Version} for all origin-symmetric convex bodies. For $p=1$, this is \cite[Lemma~3.1]{Kol_Liv_GZ}.

We begin with the observation that the local inequality takes a particularly convenient form when expressed in terms of the Riemannian structure discussed above. The following algebraic reformulation does not require symmetry.

\begin{proposition}\label{prop:transformed_local}
    Let $K\in\mathcal K^n_{2,+}$ contain the origin in its interior. For a given $f\in C^2(\partial K)$, set $\rho=f/\hbr_K$. Then \eqref{eq:local (p,q)-Version}, evaluated at $f$, can be further rewritten as
    \begin{equation}\label{eq:transfored local (p,q)}
    \begin{aligned}
        &\int_{\partial K}
            (-\Delta_{K,\mu}\rho)\rho\di V_{K,\mu}
        +\int_{\partial K}
            \big[\langle\nabla u,x\rangle-(n-p)\big]
            \rho^2\di V_{K,\mu}+(n-q)
            \left(\int_{\partial K}\rho\di V_{K,\mu}\right)^2
        \geq0.
    \end{aligned}
    \end{equation}
    Equivalently,
    \begin{equation}
    \begin{aligned}
        &\int_{\partial K}
            \left[-\Delta_{K,\mu}\rho
                +\langle\nabla u,x\rangle\rho\right]
            \rho\di V_{K,\mu} \geq (n-p)\mathrm{Var}_{V_{K,\mu}}(\rho)-(p-q)\left(\int_{\partial K}\rho\di V_{K,\mu}\right)^2,
    \end{aligned}
    \end{equation}
    where, here and below with a slight abuse of notation,
    \[
        \mathrm{Var}_{V_{K,\mu}}(v)
        \coloneqq\int_{\partial K}v^2\di V_{K,\mu}
        -\left(\int_{\partial K}v\di V_{K,\mu}\right)^2.
    \]
\end{proposition}

\begin{proof}

First, consider the substitution $f=\rho \hbr$. Observe that, by the divergence theorem, 
\[ \int_{\partial K} \langle {\rm I\! I}^{-1} \nabla_{\partial K} f, \nabla_{\partial K}f \rangle \di \mu=-\int_{\partial K} \nabla_{\partial K}\cdot ({\rm I\!I}^{-1} \nabla_{\partial K} f) f \di \mu+\int_{\partial K}  \langle {\rm I\!I}^{-1} \nabla_{\partial K} f, \nabla_{\partial K} u \rangle f \di \mu. \] 
Now, by direct computations, 
    \begin{align*}
        \nabla_{\partial K}\cdot ({\rm I\!I}^{-1} \nabla_{\partial K} f) = & \nabla_{\partial K}\cdot (\hbr_K{\rm I\!  I}^{-1} \nabla_{\partial K} \rho)+\nabla_{\partial K} \cdot (\rho {\rm I\! I}^{-1} \nabla_{\partial K} \hbr_K )\\
        = & \nabla_{\partial K} \cdot ({\rm I\! I}^{-1} \nabla_K \rho) \hbr_K + 2 \langle {\rm I\! I}^{-1} \nabla_{\partial K}\rho, \nabla_{\partial K} \hbr_K \rangle + \rho(n-1 - {\rm tr}({\rm I\! I})),
    \end{align*}
    and 
    \begin{align*}
        \langle {\rm I\!I}^{-1} \nabla_{\partial K} f, \nabla_{\partial K} u \rangle=\langle {\rm I\!I}^{-1} \nabla_{\partial K} \hbr_K, \nabla_{\partial K} u \rangle \rho + \langle {\rm I\!I}^{-1} \nabla_{\partial K} \hbr_K, \nabla_{\partial K} u \rangle \hbr_K.
    \end{align*}
    Recalling that 
    \[{\rm H}_\mu={\rm tr}({\rm I\! I})-\langle \nabla u, \nu_K \rangle\] 
    and 
    \[ \langle \nabla u, x \rangle=\langle \nabla u, \nu_K \rangle \hbr_K + \langle \nabla_{\partial K} u , {\rm I\! I}^{-1}\nabla_{\partial K} \hbr_K \rangle,  \] 
    we see that the generalized mean curvature cancels and we are left with precisely with \eqref{eq:transfored local (p,q)}. 
\end{proof}\noindent

The log-concavity of $\mu$ corresponds to $p=1$ and $q=0$,
with the latter understood in the limiting sense. Since
log-concavity holds without symmetry assumptions, this gives
the following reformulation of the boundary Poincar\'e
inequalities of Colesanti~\cite{Colesanti_poincare} and
Kolesnikov--Milman~\cite{Kol_Mil_BL,Kol_Mil_Po}. An equivalent result was independently obtained in \cite[Lemma 8.2]{cor-er_liv_rot}.
\begin{theorem}
    Let $\mu$ be a log-concave measure on $\R^n$ with density
    $\di\mu=e^{-u}\di x$, where $u\in C^2(\R^n)$ is convex.
    Then, for every $K\in\mathcal K^n_{2,+}$ containing the
    origin in its interior and every $\rho\in C^1(\partial K)$,
    \begin{equation}\label{eq:Centroaffine_poincare}
    \begin{aligned}
        &\int_{\partial K}
            \big(n-1-\langle\nabla u,x\rangle\big)
            \rho^2\di V_{K,\mu}
        -n\left(\int_{\partial K}\rho\di V_{K,\mu}\right)^2\leq
        \int_{\partial K}
            g_K(\grad_{g_K}\rho,\grad_{g_K}\rho)
            \di V_{K,\mu}.
    \end{aligned}
    \end{equation}
\end{theorem}

\begin{proof}
     It is immediate in virtue of Proposition~\ref{prop:transformed_local} with $p=1$ and $q=0$, together with \eqref{eq:IBP_w_aff}.
\end{proof}

The preceding formulations admit a natural spectral interpretation. We now return to even measures and origin-symmetric bodies and introduce the corresponding operator.

\begin{definition}[Weighted Hilbert operator]
    Let $\mu$ be an even log-concave measure on $\R^n$ with density $\di\mu=e^{-u}\di x$, where $u\in C^2(\R^n)$, and let $K\in\mathcal K^n_{2,+}$ be origin-symmetric. We define the weighted Hilbert--Brunn--Minkowski operator associated with $K$ and $\mu$ by
    \begin{equation}\label{eq: Full_Dimensional_Hilbert_operator}
        \mathcal D_{K,\mu}f
        =-\Delta_{K,\mu}f+\langle\nabla u,x\rangle f.
    \end{equation}
\end{definition} \noindent

Notice that an integration by parts gives
\[
    \int_{\partial K}
        (\mathcal D_{K,\mu}f)f\di V_{K,\mu}
    =
    \int_{\partial K}
        \left[
            g_K(\grad_{g_K}f,\grad_{g_K}f)
            +\langle\nabla u,x\rangle f^2
        \right]\di V_{K,\mu}.
\]
Therefore, $\mathcal D_{K,\mu}$ is self-adjoint in the Sobolev space $H^1(\partial K, V_{K,\mu})$. In particular, it is reasonable to introduce the following spectral quantity:
\begin{equation}\label{eq:first_even_eigenvalue}
    \lambda_{1,e}(\mathcal D_{K,\mu})
    \coloneqq \inf_{\substack{
        f\in C^2_e(\partial K)\\
        \mathrm{Var}_{V_{K,\mu}}(f)>0
    }} \mathcal{J}_{K,\mu}(f) = 
    \inf_{\substack{
        f\in C^2_e(\partial K)\\
        \mathrm{Var}_{V_{K,\mu}}(f)>0
    }}
    \frac{
        \displaystyle\int_{\partial K}
        (\mathcal D_{K,\mu}f)f\di V_{K,\mu}
    }{
        \displaystyle\mathrm{Var}_{V_{K,\mu}}(f)
    }.
\end{equation}
Here $C^2_e(\partial K)$ denotes the even $C^2$ functions on $\partial K$. Unlike in the unweighted setting \cite{Kol_Mil_Po, Mil_centroaffine}, constants do not, in general, belong to the kernel of $\mathcal D_{K,\mu}$. Consequently, one cannot impose a zero mean constraint in this variational problem. Furthermore, notice that the divergence theorem gives
\begin{equation}\label{eq:total_weighted_cone_volume}
\begin{aligned}
    V_{K,\mu}(\partial K)
    &=\frac1{n\mu(K)}
        \int_{\partial K}\hbr_K\di\mu=1-\frac1{n\mu(K)}
        \int_K\langle\nabla u,x\rangle\di\mu
    \leq1,
\end{aligned}
\end{equation}
where the last inequality follows from the convexity and evenness of $u$. In the case that $u$ is strictly convex the above inequality is strict. Indeed, since $\langle \nabla u,x \rangle \geq 0$ if $V_{K,\mu}(\partial K)=1$ we would have
\[\langle \nabla u(x)-\nabla u(0),x-0 \rangle =0, \quad \forall x \in K.\]
But then strict convexity would imply that $x=0$, which is a contradiction. Writing $m=V_{K,\mu}(\partial K)\in(0,1]$, the Cauchy--Schwarz inequality yields
\[
    \left(\int_{\partial K}f\di V_{K,\mu}\right)^2
    \leq m\int_{\partial K}f^2\di V_{K,\mu}.
\]
Thus,
\begin{equation}\label{eq:variance_lower_bound}
\begin{aligned}
    \mathrm{Var}_{V_{K,\mu}}(f)
    &=\int_{\partial K}f^2\di V_{K,\mu}
        -\left(\int_{\partial K}f\di V_{K,\mu}\right)^2 \geq(1-m)\int_{\partial K}f^2\di V_{K,\mu}\\
    &\geq\frac{1-m}{m} \left(\int_{\partial K}f\di V_{K,\mu}\right)^2 \geq0.
\end{aligned}
\end{equation}
Thus for strictly convex $u$, $\mathrm{Var}_{V_{K,\mu}}(f)>0$ for every nonzero continuous $f$. 

Finally, when one discusses concavity principles for log-concave measures it is enough to assume that the function $u$ is strictly convex, via a usual approximation argument, so we focus our attention to this case in order to avoid ambiguity.

\medskip

With these preliminaries in place, we use the above operator to set a spectral framework for studying dimensional concavity properties of even log-concave measures, as illustrated by the following characterization, which is Theorem~1.1 in the introduction.

\begin{theorem}
    Let $\mu$ be an even log-concave measure on $\R^n$ with density $\di\mu=e^{-u}\di x$, where $u\in C^2(\R^n)$. For $p\in[0,1]$, the inequality
    \begin{equation}\label{eq:dimensional_concavity_weighted}
        \mu((1-t)K+_p tL)^{\frac pn}
        \geq(1-t)\mu(K)^{\frac pn}
            +t\mu(L)^{\frac pn}
    \end{equation}
    holds for all origin-symmetric convex bodies $K,L$ and every $t\in[0,1]$ if and only if
    \[
        \lambda_{1,e}(\mathcal D_{C,\mu})\geq n-p
    \]
    for every origin-symmetric $C\in\mathcal K^n_{2,+}$. The case $p=0$ is understood in the limiting sense as
    \[
        \log\mu((1-t)K+_0tL)
        \geq(1-t)\log\mu(K)+t\log\mu(L).
    \]
    Furthermore given an origin-symmetric convex set $K$, we have that 
    \[t \mapsto \log \mu(e^{t}K),\]
    is concave at $t=0$, if and only if 
    \[\mathcal{J}_{K,\mu}(1) \geq n. \]
    In particular the (B) property holds for the measure $\mu$ if and only if $\mathcal{J}_{K,\mu}(1) \geq n$, for every origin symmetric convex set $K \in \mathcal{K}^{n}_{2,+}$.
\end{theorem}

\begin{proof}
    For the first part, setting $q=p$ in Proposition~\ref{prop:transformed_local},
    the local inequality becomes
    \[
        \int_{\partial C}
            (\mathcal D_{C,\mu}f)f\di V_{C,\mu}
        \geq(n-p)\mathrm{Var}_{V_{C,\mu}}(f).
    \]
    For functions with $\mathrm{Var}_{V_{C,\mu}}(f)>0$,
    this is exactly the bound on the Rayleigh quotient, which is the case for strictly convex potentials $u$. For $p<1$, the equivalence with the global inequality follows from Lemma 3.1 in \cite{hos_kol_liv}. For $p=1$, one can refer to \cite{Kol_Liv_GZ}.
    \medskip
    
For the second part, fixing a convex set $K$, consider the functional
    \[F(t) = \mu(e^{t}K), t \in \mathbb{R}.\]
    It is well known, see e.g. \cite{cordero2025concavity,rot_B}, that 
    \[F'(0) = \int_{K}n - \langle D u, x\rangle \di \mu = n\mu(K)\int_{\partial K} h_{K}(\nu_{K})\di V_{K,\mu} \]
    and
    \[F''(0) =  \mu(K)\int_{K} \langle Du ,x \rangle^2-\langle (D^2u) x,x \rangle -\langle Du,x \rangle \di \mu.\]
    Thus the inequality $(\log(F(t)))''(0) \leq 0$, is equivalent to
    \begin{equation}\label{eq:local_B_inequality}
    \int_{K} \langle (D^2u)x,x\rangle + \langle Du,x\rangle \di \mu \geq \int_{K} \langle D u,x \rangle^2 \di \mu -\frac{1}{\mu(K)}\left(\int_{K}\langle D u,x \rangle \di \mu \right)^{2}. 
    \end{equation}
    On the other hand, $\mathcal{J}_{K,\mu}(1) \geq n$ can be expanded as
    \begin{equation}\label{eq:local_B_boundary}
    \frac{1}{n\mu(K)}\int_{\partial K} \langle Du,x \rangle \langle x, \nu_{K}(x)\rangle \di \mu \geq \frac{1}{\mu(K)}\int_{\partial K} \langle x, \nu_{K}(x)\rangle \di \mu -\frac{1}{n\mu(K)^2}\left( \int_{\partial K} \langle x, \nu_{K}(x)\rangle \di \mu \right)^2.
    \end{equation}
    For the term on the left hand side, an integration by parts yields
    \begin{align*}
    \int_{\partial K} \langle Du,x \rangle \langle x, \nu_{K}(x)\rangle e^{-u}\di x &=  \int_{K} \langle Du,x \rangle \LB_{\mu}(\frac{|x|^2}{2})\di \mu + \int_{K}\langle D (\langle Du,x\rangle ), x\rangle \di \mu \\
    &= n\int_{K}\langle Du,x \rangle \di \mu - \int_{K} \langle Du,x\rangle^2 \di \mu \\
    &+\int_{K} \langle (D^2u)x,x\rangle \di \mu +\int_{K} \langle Du,x \rangle \di \mu
    \end{align*}
    For the right hand side of \eqref{eq:local_B_boundary}, we need only employ the following direct application of the divergence theorem
    \[\int_{\partial K} \langle x, \nu_{K}(x)\rangle \di \mu = n\mu(K)-\int_{K}\langle Du,x \rangle \di \mu, \]
    which yields that the right hand side of \eqref{eq:local_B_boundary} is equal to
    \[\int_{K} \langle Du,x \rangle \di \mu -\frac{1}{n\mu(K)}\left( \int_{K} \langle Du,x \rangle^2 \right).\]
    Putting everything together, $\frac{1}{n}$ cancels out and we are left precisely with \eqref{eq:local_B_inequality}, concluding the proof. 
\end{proof} 

\section{Anisotropic Poincar\'e and strong Brascamp--Lieb inequalities for gauge-symmetric measures}

We start with a simple lemma to split the Hessian of inverse quadratic gauge potentials into radial and tangential components, with respect to the anisotropic Riemannian metric introduced in \cite{Ivaki_Milman_p}. We remark that a proof can be given via a usage of Lemma~2.2 therein. See, for example, \cite[Lemma 8.1]{cor-er_liv_rot}. We propose an independent proof.

\begin{lemma}\label{lemma:centroaffine_decomposition}
For $K$ smooth convex body, consider the function $\Phi(x)=\frac{\|x\|_K^2}{2}$. For every $\varphi \in C^1(\R^n)$ and $y \in \partial K$, we have 
\[ \langle (D^2 \Phi(y))^{-1} D \varphi(y), D \varphi(y) \rangle=\langle D \varphi(y), y \rangle^2+g_K({\rm grad}_{g_K} \varphi(y), {\rm grad}_{g_K} \varphi(y)).\]
\end{lemma}
\begin{proof}
    We first observe that, by Fenchel--Legendre duality, 
    \[ \left(D^2 \Phi(x)\right)^{-1}=D^2 \left( \frac{h_K^2}{2}\right)(D \Phi(x))\] 
    for every $x \in \R^n \setminus\{o\}$. The equality above can be found, for example, in \cite[(1.49), Section 1.7]{Schneider_CB}. Then, direct computations show 
    \[ D^2 \left(\frac{h_K^2}{2}\right)(D \Phi(x))=\left.\left( h_K D^2 h_K + D h_K \otimes D h_K\right)\right|_{z=D \Phi(x)}. \] 
    If we now choose $x=y \in \partial K$, by the properties of gauge functions $D \Phi(y)=\nu_K(y)/\hbr_K(y)$. Denoting by $\pi_z, z \in \R^n\setminus\{o\},$ the orthogonal projection along $z$, it is well-known that (compare, e.g., \cite[Lemma 2.5.1]{Schneider_CB}) 
    \[ D^2 h_K(z)=|z|^{-1} \pi_z^t D^2 h_K(z/|z|) \pi_z,\] 
    where now we consider these objects as bilinear forms on $\R^n$. In particular, 
    \[ D^2 h_K(\nu_K(y)/\hbr(y))=\hbr_K(y)\left(\pi_{\nu_K(y)}^t D^2 h_K(\nu_K(y)) \pi_{\nu_K(y)}\right).\]
    As $h_K$ is $1$-homogeneous and $\nabla h_K$ is $0$-homogeneous, we infer
    \[(D^2 \Phi(y))^{-1}=\hbr_K(y)\left(\pi_{\nu_K(y)}^t D^2 h_K(\nu_K(y)) \pi_{\nu_K(y)}\right)+y \otimes y,\] 
    and the proof is concluded by direct computations, recalling that, for every $y \in \partial K$, 
    \[ D^2 h_K(\nu_K(y))={\rm I\! I}^{-1}(y).\] 
\end{proof}
\noindent
Finally, let us recall the following $1$-dimensional Poincar\'e inequality, which is \cite[Lemma 8]{rot_B}.
\begin{lemma}\label{lemma:1_po}
    Consider $v,w \in C([0,\infty))\cap C^1((0,\infty))$, and $f \in C^1_c([0,\infty))$ such that $f(0)=0$. Then, for every $\alpha\geq 0$ 
    \[ \int_0^\infty \frac{w'}{t}f^2 t^\alpha e^{-w-v}\di t \leq \int_0^\infty \left((f')^2 + \alpha\left(\frac{f}{t}\right)^2-v'\frac{f^2}{t} \right)t^\alpha e^{-w-v}\di t.\]
\end{lemma}

With these results at hand, we can now prove the main result of this section, which can be understood as the anisotropic analogue of \cite[Theorem 4]{rot_B}. An equivalent formulation was recently proved in \cite[Theorem 1.8]{cor-er_liv_rot}.
\begin{theorem}\label{thm:odd_BL}
    Consider an increasing function $w \in C^1(0,\infty)$ and $K$ an even and smooth convex body. Let $u: \mathbb{R}^n \rightarrow \mathbb{R}$ be convex, then denoting by $\eta$ the measure on $\R^n$ satisfying $\di \eta(x)=e^{-u(x)-w(\|x\|_K)}\di x$, for every odd function $\varphi \in C^1(\R^n)$ we have 
    \begin{equation}\label{eq:an_po}
        \int_{\R^n} \frac{w'(\|x\|_K)}{\|x\|_K} \varphi^2(x) \di \eta(x) \leq \int_{\R^n} \left \langle \left( D^2 \Phi(x) \right)^{-1} D \varphi(x), D \varphi(x) \right \rangle \di \eta(x),
    \end{equation}
    where $\Phi(x)=\|x\|_K^2/2$.
\end{theorem}
\begin{proof}
    We start by observing that, using polar coordinates, we can rewrite the left-hand side of \eqref{eq:an_po} as
    \begin{equation}\label{eq:polar_coordinates}
        \int_{\R^n} \frac{w'(\|x\|_K)}{\|x\|_K} \varphi^2(x) \di \eta(x)=\int_{\partial K} \int_0^\infty \frac{w'(r)}{r} \varphi^2(ry)r^{n-1}e^{-u(ry)-w(r)} \langle \nu_K(y),y \rangle\di r \di \haus^{n-1}(y).
    \end{equation}
    By Lemma \ref{lemma:1_po} with the choice $\alpha=n-1$, we infer
    \begin{align} \notag
    \text{\eqref{eq:polar_coordinates}}\leq  \int_{\partial K} \int_0^\infty &\left( \langle D \varphi(ry),y\rangle^2+(n-1)\left( \frac{\varphi(ry)}{r}\right)^2 - \frac{1}{r}\langle Du(ry),ry\rangle\frac{\varphi(ry)^2}{r} \right)\\
    &\times r^{n-1} e^{-w(r)-u(ry)}\langle \nu_K(y),y \rangle \di r \di \haus^{n-1}(y). \label{eq:red_1_an}
    \end{align}
    Next, observe that, by \eqref{eq:Centroaffine_poincare}, since $h$ is odd and the measure as well as $\partial K$ is symmetric, we have that 
    \begin{align} \notag
        \int_{\partial K}  \left( (n-1)\left( \frac{\varphi(ry)}{r}\right)^2 - \langle Du(ry),ry\rangle\left( \frac{\varphi(ry)}{r}\right)^2\right)\langle \nu_K(y),y \rangle e^{-u(ry)} &\di \haus^{n-1}(y)\\
        \leq  \int_{\partial K}  g_K\left({\rm grad}_{g_K} \left( \frac{\varphi(ry)}{r} \right), {\rm grad}_{g_K} \left( \frac{\varphi(ry)}{r} \right) \right) \langle \nu_K(y),y \rangle e^{-u(ry)} &\di \haus^{n-1}(y). \label{eq:red_2_an}
    \end{align}
    If now we invert the order of integration in \eqref{eq:red_1_an} and apply the inequality \eqref{eq:red_2_an}, we deduce
    \begin{align*}
        \text{\eqref{eq:red_1_an}}\leq & \int_{\partial K}\int_0^\infty \left(  \langle D \varphi(ry),y\rangle^2+ \right. \\ &\left. g_K\left({\rm grad}_{g_K} \left( \frac{\varphi(ry)}{r} \right), {\rm grad}_{g_K} \left( \frac{\varphi(ry)}{r} \right) \right) \right)r^{n-1}\langle \nu_K(y),y \rangle e^{-u(ry)-w(r)}\di r \di \haus^{n-1}(y).
    \end{align*}
    Finally, by Lemma \ref{lemma:centroaffine_decomposition} we have 
    \[ \langle D \varphi(ry),y\rangle^2+ g_K\left({\rm grad}_{g_K} \left( \frac{\varphi(ry)}{r} \right), {\rm grad}_{g_K} \left( \frac{\varphi(ry)}{r} \right) \right) = \left \langle \left( D^2 \Phi(ry) \right)^{-1} D \varphi(ry), D \varphi(ry) \right \rangle, \] 
    and the proof is concluded by integrating in polar coordinates.  
\end{proof}

\begin{remark}
    The above inequality is sharp whenever $V=0$ and $h=\langle \nabla \Phi, \theta \rangle$, with $\theta$ fixed. We remark also that the above inequality is a Poincar\'e inequality on the Hessian weighted manifold
    \[(\mathbb{R}^n, D^2 \Phi, e^{-V-w(\|x\|_{K})}\di x)\]
    for odd functions. In the case that the underlying set $K$ is a ball, we have that $D^2\Phi$ is the identity and thus the above inequality reduces to a Poincar\'e inequality in the usual Euclidean space, which is the result of Cordero-Erausquin and Rotem \cite{rot_B}.
\end{remark}

The usage of the centroaffine metric (or the ambient anisotropic one) albeit convenient for the computations and the proof of the inequality still does not agree fully with the geometry of $\mathbb{R}^n$. This problem is remedied in the works of Milman and Ivaki--Milman \cite{Mil_centroaffine,Ivaki_Milman_p}, through a pinching condition on the set $K$. Following  this idea we establish an improved Brascamp--Lieb inequality. The proof is an adaptation of \cite[Theorem 3]{rot_B}.
\begin{theorem}\label{thm:strong_BL}
    Consider an increasing function $w \in C^2((0,\infty))$ and a probability measure $\eta$ on $\R^n$ satisfying $\di \eta(x)=e^{-V(x)-w(\|x\|_K)}\di x$, where $V \in C^2(\R^n)$ is convex and even and $K$ is an origin-symmetric smooth convex body. Suppose, additionally, that \[ \alpha \Id \leq D^2 \frac{\|x\|_K^2}{2} \leq \beta \Id\] for every $x \in \R^n$ some $0<\alpha<\beta<\infty$, and that $w$ satisfies
    \begin{equation}\label{eq:cond_pos_BL}
         w''(s)+\gamma \frac{w'(s)}{s} \geq 0, \quad s >0
    \end{equation}
    where $\gamma=\alpha/\beta \leq 1$.\\
    Then, writing $W(x)=w(\|x\|_K)$, for every even function $\varphi \in C^1(\R^n)$ we have the inequality 
    \begin{equation}\label{eq:strong_BL}
        {\rm Var}_\eta(\varphi) \leq \int_{\R^n} \left \langle \left(D^2 W(x)+\gamma \frac{w'(\|x\|_K)}{\|x\|_K}D^2 \frac{\|x\|_K^2}{2} \right)^{-1} D \varphi(x), D \varphi(x)  \right \rangle \di \eta(x).
    \end{equation}
\end{theorem}
\begin{proof}
    The arguments mimic those in \cite{rot_B} which, in turn, builds upon \cite{Cor_Fra_Mau}. We streamline the presentation for the convenience of the reader. 
    
    Let us begin by showing that \eqref{eq:cond_pos_BL} implies that 
    \[ \text{the matrix }\A(x)=D^2 W(x)+\gamma \frac{w'(\|x\|_K)}{\|x\|_K}D^2 \frac{\|x\|_K^2}{2} \text{ is positive semidefinite.}\] 
    Indeed, by direct computations, for every $x\in \R^n\setminus \{o\}$ and writing $t=\|x\|_K$, we have the identity 
    \begin{align*}
       \A(x)=& w'(t)D^2\|x\|_K+w''(t) D\|x\|_K \otimes D\|x\|_K+\gamma  \frac{w'(t)}{t} \left(tD^2\|x\|_K+D\|x\|_K \otimes D\|x\|_K \right)\\
        =& \left(w''(t)+ \gamma  \frac{w'(t)}{t}\right)D\|x\|_K \otimes D\|x\|_K+w'(t)(1+\gamma)D^2\|x\|_K.
    \end{align*}
    The coefficient of the first matrix is non-negative thanks to \eqref{eq:cond_pos_BL} while, for the second one, this follows from the monotonicity of $w$ and, thus, $\B$ is positive semidefinite. This is enough to follow through with the computations in the rest of the proof, as explained in \cite[Remark 9]{rot_B}. We may (and do) assume that $\varphi$ is in addition centered, i.e. $\int \varphi=0$ subtracting a constant if that is not the case. With these preliminary remarks in place notice that \eqref{eq:strong_BL} is equivalent to
    \begin{equation}\label{eq:red_1_BL}
        \int_{\R^n} \langle \A^{-1} D\varphi, D \varphi \rangle-\varphi^2 \di \eta \geq 0.
    \end{equation}
    Denoting by $\LB_\eta$ the standard Laplace--Beltrami operator on $\R^n$ with respect to $\eta$, $\LB_\eta(C^2_c(\R^n))$ is dense in the subspace of $L^2(\eta)$ orthogonal to constant functions. See \cite[Lemma 3]{Cor_Fra_Mau} for a proof. In particular, we can approximate $\varphi$ with these function and, therefore, \eqref{eq:red_1_BL} would follow if we are able to show that
    \begin{align} \notag
        &\int_{\R^n}  \langle \A^{-1} D\varphi, D \varphi \rangle-2\LB_\eta(u)\varphi+(\LB_\eta(u))^2 \di \eta \\
        =&\int_{\R^n} \langle \A^{-1} D\varphi, D \varphi \rangle-\varphi^2 + (\LB_\eta(u)-\varphi)^2 \di \eta \geq 0 \quad \forall u \in C^2_c(\R^n). \label{eq:red_2_BL}
    \end{align}
    
    We now estimate the first integral in \eqref{eq:red_2_BL}. We have $-\int_{\R^n} \LB_\eta(u)\varphi \di \eta=\int_{\R^n} \langle D u, D \varphi \rangle \di \eta$ by integration by parts. The Cauchy--Schwarz inequality for positive quadratic forms implies that
    \[\langle \A^{-1} D\varphi, D \varphi \rangle+\langle \A Du, Du \rangle \geq -2\langle Du, D\varphi \rangle,\] we obtain the bound
    \begin{equation}\label{eq:red_3_BL}
        \int_{\R^n}  \langle \A^{-1} D\varphi, D \varphi \rangle-2\LB_\eta(u)\varphi \di \eta \geq \int_{\R^n} -  \langle \A Du, Du \rangle \di \eta.
    \end{equation}
    Concerning the remainder, integration by parts gives
    \begin{equation}\label{eq:red_4_BL}
         \int_{\R^n} (\LB_\eta(u))^2 \di \eta = \int_{\R^n} \| D^2 u \|^2_{\rm HS}+\langle (D^2 V+D^2 W) Du, Du\rangle \di \eta,
    \end{equation}
    where $\| \cdot \|_{\rm HS}$ denotes the Hilbert--Schmidt norm of a matrix. Then, \eqref{eq:red_3_BL} and \eqref{eq:red_4_BL} in combination with the fact that $D^2V$ is positive definite, give 
    \begin{align*}
        \int_{\R^n}  \langle \A^{-1}  D\varphi, D \varphi \rangle&-2\LB_\eta(u)\varphi+(\LB_\eta(u))^2 \di \eta \geq \\ &\int_{\R^n} \| D^2 u \|^2_{\rm HS}-\left \langle \gamma \frac{w'(\|x\|_K)}{\|x\|_K}D^2 \frac{\|x\|_K^2}{2} D u, D u \right  \rangle \di \eta.
    \end{align*}
    Therefore to conclude the proof, it suffices to show that 
    \[ \int_{\R^n} \| D^2 u \|^2_{\rm HS} \di \eta \geq  \int_{\R^n} \gamma \frac{w'(\|x\|_K)}{\|x\|_K} \left \langle D^2 \frac{\|x\|_K^2}{2} D u, D u \right  \rangle \di \eta \quad \forall u \in C^2_c(\R^n), \] which is indeed true, as
    \begin{align}\label{eq:spectral_estimate_Hessian}
        \int_{\R^n} \| D^2 u \|^2_{\rm HS} \di \eta \geq &\alpha \int_{\R^n}  \sum_{i=1}^n \left \langle  \left(D^2 \frac{\|x\|_K^2}{2}\right)^{-1} D u_i, D u_i  \right \rangle \di \eta \\ \geq & \alpha\int_{\R^n}  \frac{w'(\|x\|_K)}{\|x\|_K} |D u |^2 \di \eta  \geq \int_{\R^n} \gamma  \frac{w'(\|x\|_K)}{\|x\|_K} \left \langle D^2 \frac{\|x\|_K^2}{2} D u, D u \right  \rangle \di \eta,
    \end{align}
    where $u_i$ denotes the $i$-th derivative of $u$ and, in the second step, we have used Theorem \ref{thm:odd_BL}.
\end{proof}

We also have the following variant of the above inequality that can be deduced by a direct modification of the final step in the proof above \eqref{eq:spectral_estimate_Hessian}, when the function $w$ is assumed to be in addition convex. 

\begin{theorem}\label{thm:strong_BL_convex}
    Consider an increasing and convex function $w \in C^2((0,\infty))$ and a probability measure $\eta$ on $\R^n$ satisfying $\di \eta(x)=e^{-V(x)-w(\|x\|_K)}\di x$, where $V \in C^2(\R^n)$ is convex and even and $K$ is an origin-symmetric smooth convex body and set $W(x)=w(\|x\|_K)$. Assuming that
    \[D^{2} \frac{\|x\|^2_{K}}{2} \geq \alpha \mathrm{Id},\]
    then for every even function $\varphi \in C^{1}(\mathbb{R}^n)$,
        \begin{equation}\label{eq:strong_BL_v2}
        {\rm Var}_\eta(\varphi) \leq \int_{\R^n} \left \langle \left(D^2 W(x)+\alpha \frac{w'(\|x\|_K)}{\|x\|_K}\mathrm{Id} \right)^{-1} D \varphi(x), D \varphi(x)  \right \rangle \di \eta(x).
    \end{equation}
    
\end{theorem}

\begin{remark}
    Whenever $K$ is a $C^2$ smooth convex body that is origin symmetric, with strictly positive principal curvatures we may always find numbers $\alpha,\beta$ such that the above pinching or one sided control is satisfied.
\end{remark}

\section{Improved dimensional concavity for gauge-symmetric and pinched measures}

The first aim of this section is to prove a dimensional Brunn--Minkowski inequality for the following type of measures. We consider $W(x)=w(\|x\|_K)$, where $w$ is chosen as in Theorem \ref{thm:strong_BL}, and the corresponding measure $\mu$ such that $\di\mu(x)=e^{-W(x)}\di x$. Again with the notation $\Phi(x)=\|x\|^2_K/2$, by homogeneity $\Phi$ we have $D\Phi(x)=(D^2 \Phi(x))x$ and the Laplace--Beltrami operator of $\mu$ can be written as 
\[ \LB_\mu(\varphi)=\Delta \varphi - \frac{w'(\|x\|_K)}{\|x\|_K}\left \langle (D^2 \Phi)x, D\varphi \right\rangle.\] 
Consider now the measure $\eta$ in Theorem \ref{thm:strong_BL} and define, for $\varphi, \psi \in C_c^2(\R^n)$, the bilinear form 
\[ Q_\eta(\varphi,\psi)\coloneqq \int_{\R^n} \tra\left( D^2 \varphi \left( D^2 \Phi\right)^{-1} D^2 \psi\right) - \frac{w'(\|x\|_K)}{\|x\|_K}\langle D \varphi, D \psi \rangle \di \eta.\] 
Observe that, by Theorem \ref{thm:odd_BL}, $Q_\eta$ is positive semidefinite on the subspace of even functions. The operators $Q_\eta$ and $\LB_\mu$ are connected by the relation 
\begin{equation}\label{eq:bil_LB}
    Q_\eta(\varphi,\Phi)=\int_{\R^n} \LB_\mu(\varphi) \di \eta,
\end{equation}
which is instrumental in the proof of the following result.
\begin{theorem}\label{thm:her_conv}
    Consider an origin-symmetric and smooth convex body $K\subset \R^n$ such that the corresponding quadratic potential $\Phi$ satisfies
    \begin{equation}\label{eq:pinching}
        \alpha \Id \leq D^2 \Phi \leq \beta \Id
    \end{equation}
    for some $0 < \alpha \leq \beta <\infty$. Consider, moreover, an increasing function $w \in C^2((0,\infty))$ satisfying \eqref{eq:cond_pos_BL}, and write $W(x)=w(\|x\|_K)$. Then, for every convex function $V \in C^2(\R^n)$ such that the measure $\eta$ with density $\di \eta(x)=e^{-V(x)-W(x)}\di x$ is a probability measure we have 
    \begin{equation}\label{eq:gamma_her_conv}
        \int_{\R^n} \|D^2 \varphi \|_{\rm HS}^2 + \left\langle (D^2 W) D \varphi, D \varphi\right \rangle \di \eta \geq \frac{\gamma}{n} \left( \int_{\R^n} \LB_\mu(\varphi) \di \eta\right)^2
    \end{equation}
    for every even $\varphi \in C^2_c(\R^n)$.
\end{theorem}
\begin{proof}
    Fix $\varphi \in C^2_c(\R^n)$ even. We start by observing that, by the pinching condition \eqref{eq:pinching}, 
    \begin{align} \notag
        &\int_{\R^n} \| D^2 \varphi \|^2_{\rm HS}+\left\langle (D^2 W) D \varphi, D \varphi\right \rangle \di \eta \geq \int_{\R^n}\alpha \tra(D^2 \varphi (D^2 \Phi)^{-1} D^2 \varphi )+\left\langle (D^2 W) D \varphi, D \varphi\right \rangle \di \eta\\
         & \geq  \alpha Q_\eta(\varphi,\varphi)+ \int_{\R^n} \left \langle \left( D^2 W + \gamma \frac{w'(\|x\|_K)}{\|x\|_K}D^2 \Phi \right) D \varphi, D \varphi \right \rangle \di \eta \geq \alpha Q_\eta(\varphi,\varphi),\label{eq:red_1_her_conv}
    \end{align}
    where the last step follows from \eqref{eq:cond_pos_BL} as in Theorem \ref{thm:strong_BL}.

    Next, we note that, by the bilinearity of $Q_\eta$ and \eqref{eq:bil_LB}, the function $\psi=\varphi+\lambda \Phi, \lambda \in \R,$ satisfies 
    \[ Q_\eta(\psi,\psi)=Q_\eta(\varphi,\varphi)+2\lambda \int_{\R^n} \LB_\mu(\varphi) \di \eta + \lambda^2 \int_{\R^n} \LB_\mu(\Phi)\di \eta.\] 
    Then, the assumed estimate $\LB_\mu(\Phi) \leq  \beta n$ implies 
    \[ Q_\eta(\varphi,\varphi) \geq Q_\eta(\psi,\psi)-2\lambda \int_{\R^n} \LB_\mu(\varphi) \di \eta -\lambda^2  \beta n \] for every $\lambda \in \R$. Optimizing in $\lambda$ yields 
    \[ Q_\eta(\varphi,\varphi) \geq Q_\eta(\psi,\psi)+\frac{1}{\beta n} \left( \int_{\R^n} \LB_\mu(\varphi)\di \eta\right)^2 \geq \frac{1}{\beta n} \left( \int_{\R^n} \LB_\mu(\varphi)\di \eta\right)^2, \] 
    where the last line follows since $Q_\eta$ is positive semidefinite for the functions at hand. The proof of \eqref{eq:gamma_her_conv} is concluded comparing this last estimate with \eqref{eq:red_1_her_conv}.
\end{proof}

To show that a dimensional Brunn--Minkowski inequality holds for the measure $\mu$ with concavity $\gamma/n$, we invoke the following result by Kolesnikov and Livshyts, which is part of \cite[Lemma 2.3]{Kol_Liv_GZ}.
\begin{theorem}\label{thm:her_implies_dim}
    Let $\overline{\mu}$ be an even measure with density $\di \mu(x)=e^{-U(x)}\di x$ such that $U \in C^2(\R^n)$, suppose that for every origin symmetric convex body $A \subset \R^n$ and every $\varphi \in C^2(\R^n)$ with $\LB_{\overline{\mu}}(\varphi)=1$ on $A$ we have that
    \begin{equation}\label{eq:conc_power_her}
        \int_{A} \|D^2 \varphi \|_{\rm HS}^2+ \langle (D^2 U)D \varphi, D \varphi \rangle \di \overline{\mu} \geq p.
    \end{equation}
    Then,
    \[ \mu((1-t)A+tB)^{p} \geq (1-t)\mu(A)^{p}+t\mu(A)^{p},\]
    for all origin-symmetric convex bodies $A,B$.
\end{theorem}\noindent
With Theorem \ref{thm:her_implies_dim} at hand, the inequality
\[ \mu((1-t) A + t B)^{\gamma/n} \geq (1-t)\mu(A)^{\gamma/n}+t \mu(B)^{\gamma/n}\] 
for all origin-symmetric convex bodies $A,B \subset \R^n$ follows by \eqref{eq:gamma_her_conv} and standard approximation arguments, providing a proof of Theorem \ref{thm:intro_dimensional_gauge}.
\medskip

We finish this section by exporting the perturbative argument in the proof of Theorem \ref{thm:her_conv} to more general pinched log-concave measures. Its proof is completely independent from most of the machinery developed above, but it came as a natural consequence of the procedures in the previous sections. The result improves the dependence on the pinching as established in \cite{Kol_Liv_GZ} and \cite{Armen_Uli}.
\begin{theorem}
    Consider an even and convex function $V \in C^2(\R^n)$ such that \[ k_1 \Id \leq D^2 V \text{ and } \Delta V \leq k_2 n\] for some $0<k_1 \leq k_2 <\infty$, and write $r=k_1/k_2$. Denoting by $\overline{\mu}$ the log-concave measure with density $\di \overline{\mu}(x)=e^{-V(x)}\di x$ we have, for every $t \in [0,1]$,
    \begin{equation}\label{eq:dim_pinched}
         \overline{\mu}((1-t)A+tB)^\frac{r}{n} \geq (1-t) \overline{\mu}(K)^{\frac{r}{n}} + t \overline{\mu}(L)^{\frac{r}{n}}
    \end{equation}
    for every origin-symmetric convex bodies $A,B \subset \R^n$.
\end{theorem}
\begin{proof}
    In order to prove \eqref{eq:dim_pinched}, the strategy is again to prove that condition \eqref{eq:conc_power_her} holds for every origin-symmetric and smooth convex body $A \subset \R^n$ and apply the conclusion of Theorem \ref{thm:her_implies_dim}. The statement follows then by density.

    Fix a smooth and origin symmetric convex body $A$. Similarly as we did with $Q_\eta$ at the beginning of this section, for $\varphi, \psi \in C^2_c(\R^n)$ we define the bilinear form 
    \[ R_V(\varphi, \psi)=\int_A \tra \left( D^2 \varphi \left( D^2 V \right)^{-1} D^2 \psi \right) - \langle D \varphi, D \psi \rangle \di \overline{\mu}.\] 
    Observe that $R_V$ is positive definite on odd functions by the Brascamp--Lieb variance inequality. For $\varphi$ fixed, consider the linear perturbation $\psi=\varphi + \lambda V, \lambda \in \R$. Then, 
    \[ R_V(\psi, \psi)= R_V(\varphi, \varphi)+2 \lambda \int_{A} \LB_{\overline{\mu}}(\varphi) \di \overline{\mu} + \lambda^2 \int_{A} \LB_{\overline{\mu}}(V) \di \overline{\mu} \] 
    and, since $\Delta V \leq k_2 n$, we infer 
    \[ R_V(\varphi,\varphi) \geq R_V(\psi, \psi) -2 \lambda \int_{A} \LB_{\overline{\mu}}(\varphi) \di \overline{\mu} - \lambda^2 k_2 n\]  for every $\lambda \in \R$. Optimizing in $\lambda$, 
    \begin{equation}\label{eq:red_1_dim}
        R_V(\varphi, \varphi) \geq R_V(\psi, \psi) + \frac{1}{k_2 n}\left( \int_{A} \LB_{\overline{\mu}}(\varphi)\di \overline{\mu} \right)^2 \geq  \frac{1}{n k_2}\left( \int_{A} \LB_{\overline{\mu}}(\varphi)\di \overline{\mu} \right)^2.
    \end{equation}

    Finally, using $k_1 \Id \leq D^2 V$ and \eqref{eq:red_1_dim}, we estimate
    \begin{align*}
        & \int_{A} \| D^2 \varphi \|_{\rm HS}^2 +\langle (D^2 V) D\varphi, D \varphi \rangle \di \overline{\mu} \geq  \int_{A} \| D^2 \varphi \|_{\rm HS}^2 \di \overline{\mu} \\
         \geq & k_1\int_{A} \tra \left( D^2 \varphi \left( D^2 V \right)^{-1} D^2 \varphi \right) \geq k_1 R_V (\varphi,\varphi) \geq  \frac{r}{n}\left( \int_{A} \LB_{\overline{\mu}}(\varphi)\di \overline{\mu} \right)^2,
    \end{align*}
    concluding the proof.
\end{proof}

\medskip

\paragraph{AI use disclosure.} The authors used AI to streamline the organization of the paper and, in part, to improve the presentation. All the mathematical statements and proofs have been obtained and checked solely by the authors, without the aid of LLMs.

\paragraph{Acknowledgments.} The first-named author was supported by the Austrian Science Fund (FWF): 10.55776/P36344N. The second-named author gratefully acknowledges the support of the Basic Research Program of HSE University. The third-named author was supported, in part, by the Austrian Science Fund (FWF): 10.55776/PAT3787224. The authors extend their gratitude to the Erd\H{o}s Center in Budapest, to its director, K\'aroly B\"or\"oczky, and to the administrative staff for hosting and supporting the ``Focused Workshop on the log-Brunn-Minkowski conjecture: State of the art and ramifications'', during which the authors started the collaboration leading to this manuscript.

\footnotesize
\bibliography{references}

@article {Borell,
    AUTHOR = {Borell, Christer},
     TITLE = {Convex set functions in {$d$}-space},
   JOURNAL = {Period. Math. Hungar.},
  FJOURNAL = {Periodica Mathematica Hungarica. Journal of the J\'anos Bolyai
              Mathematical Society},
    VOLUME = {6},
      YEAR = {1975},
    NUMBER = {2},
     PAGES = {111--136},
      ISSN = {0031-5303,1588-2829},
   MRCLASS = {26A86 (60A10)},
  MRNUMBER = {404559},
MRREVIEWER = {M.\ D.\ Perlman},
       DOI = {10.1007/BF02018814},
       URL = {https://doi.org/10.1007/BF02018814},
}

@article {Colesanti_poincare,
    AUTHOR = {Colesanti, Andrea},
     TITLE = {From the {B}runn--{M}inkowski inequality to a class of
              {P}oincar\'e-type inequalities},
   JOURNAL = {Commun. Contemp. Math.},
  FJOURNAL = {Communications in Contemporary Mathematics},
    VOLUME = {10},
      YEAR = {2008},
    NUMBER = {5},
     PAGES = {765--772},
      ISSN = {0219-1997,1793-6683},
   MRCLASS = {52A20 (26D10)},
  MRNUMBER = {2446898},
       DOI = {10.1142/S0219199708002971},
       URL = {https://doi.org/10.1142/S0219199708002971},
}

@article {rot_B,
    AUTHOR = {Cordero-Erausquin, Dario and Rotem, Liran},
     TITLE = {Improved log-concavity for rotationally invariant measures of
              symmetric convex sets},
   JOURNAL = {Ann. Probab.},
  FJOURNAL = {The Annals of Probability},
    VOLUME = {51},
      YEAR = {2023},
    NUMBER = {3},
     PAGES = {987--1003},
      ISSN = {0091-1798,2168-894X},
   MRCLASS = {52A40 (26D10 60E15)},
  MRNUMBER = {4583060},
MRREVIEWER = {Tomasz\ Tkocz},
       DOI = {10.1214/22-aop1604},
       URL = {https://doi.org/10.1214/22-aop1604},
}

@article {Kol_Mil_BL,
    AUTHOR = {Kolesnikov, Alexander V. and Milman, Emanuel},
     TITLE = {Brascamp-{L}ieb-type inequalities on weighted {R}iemannian
              manifolds with boundary},
   JOURNAL = {J. Geom. Anal.},
  FJOURNAL = {Journal of Geometric Analysis},
    VOLUME = {27},
      YEAR = {2017},
    NUMBER = {2},
     PAGES = {1680--1702},
      ISSN = {1050-6926,1559-002X},
   MRCLASS = {53C21 (58J32 58J50)},
  MRNUMBER = {3625169},
MRREVIEWER = {Mohammed\ El A\"idi, Universidad Nacional de Colombia},
       DOI = {10.1007/s12220-016-9736-5},
       URL = {https://doi.org/10.1007/s12220-016-9736-5},
}

@article {Kol_Mil_p,
    AUTHOR = {Kolesnikov, Alexander V. and Milman, Emanuel},
     TITLE = {Local {$L^p$}-{B}runn-{M}inkowski inequalities for {$p<1$}},
   JOURNAL = {Mem. Amer. Math. Soc.},
  FJOURNAL = {Memoirs of the American Mathematical Society},
    VOLUME = {277},
      YEAR = {2022},
    NUMBER = {1360},
     PAGES = {v+78},
      ISSN = {0065-9266,1947-6221},
      ISBN = {978-1-4704-5160-8; 978-1-4704-7092-0},
   MRCLASS = {52A40 (35P15 52A23 58J50)},
  MRNUMBER = {4438690},
MRREVIEWER = {Arnaud\ Marsiglietti},
       DOI = {10.1090/memo/1360},
       URL = {https://doi.org/10.1090/memo/1360},
}

@article {Kol_Mil_Po,
    AUTHOR = {Kolesnikov, Alexander V. and Milman, Emanuel},
     TITLE = {Poincar\'e and {B}runn--{M}inkowski inequalities on the
              boundary of weighted {R}iemannian manifolds},
   JOURNAL = {Amer. J. Math.},
  FJOURNAL = {American Journal of Mathematics},
    VOLUME = {140},
      YEAR = {2018},
    NUMBER = {5},
     PAGES = {1147--1185},
      ISSN = {0002-9327,1080-6377},
   MRCLASS = {53C21 (52A40 58D25)},
  MRNUMBER = {3862061},
MRREVIEWER = {Yunwei\ Xia},
       DOI = {10.1353/ajm.2018.0027},
       URL = {https://doi.org/10.1353/ajm.2018.0027},
}

@article{cordero2025concavity,
  title={Concavity principles for weighted marginals},
  author={Cordero-Erausquin, Dario and Eskenazis, Alexandros},
  journal={arXiv:2506.16941},
  year={2025}
}

@article{Liv_Dimensional_BM,
  title={A universal bound in the dimensional {B}runn--{M}inkowski inequality for log-concave measures},
  author={Livshyts, Galyna V},
  journal={Trans. Amer. Math. Soc.},
  volume={376},
  number={09},
  pages={6663--6680},
  year={2023}
}

@article {Sar,
    AUTHOR = {Saroglou, Christos},
     TITLE = {More on logarithmic sums of convex bodies},
   JOURNAL = {Mathematika},
  FJOURNAL = {Mathematika. A Journal of Pure and Applied Mathematics},
    VOLUME = {62},
      YEAR = {2016},
    NUMBER = {3},
     PAGES = {818--841},
      ISSN = {0025-5793,2041-7942},
   MRCLASS = {52A40 (52A20 52A39)},
  MRNUMBER = {3521355},
MRREVIEWER = {Gaoyong\ Zhang},
       DOI = {10.1112/S0025579316000061},
       URL = {https://doi.org/10.1112/S0025579316000061},
}

@article{nayar2013note,
  title={A note on a {B}runn--{M}inkowski inequality for the {G}aussian measure},
  author={Nayar, Piotr and Tkocz, Tomasz},
  journal={Proceedings of the American Mathematical Society},
  pages={4027--4030},
  year={2013},
  publisher={JSTOR}
}

@article{gardner2010gaussian,
  title={Gaussian {B}runn--{M}inkowski inequalities},
  author={Gardner, Richard and Zvavitch, Artem},
  journal={Trans. Amer. Math. Soc.},
  volume={362},
  number={10},
  pages={5333--5353},
  year={2010}
}

@article {ColLivMar_stability,
    AUTHOR = {Colesanti, Andrea and Livshyts, Galyna V. and Marsiglietti,
              Arnaud},
     TITLE = {On the stability of {B}runn-{M}inkowski type inequalities},
   JOURNAL = {J. Funct. Anal.},
  FJOURNAL = {Journal of Functional Analysis},
    VOLUME = {273},
      YEAR = {2017},
    NUMBER = {3},
     PAGES = {1120--1139},
      ISSN = {0022-1236,1096-0783},
   MRCLASS = {52A40},
  MRNUMBER = {3653949},
MRREVIEWER = {Stefano\ Campi},
       DOI = {10.1016/j.jfa.2017.04.008},
       URL = {https://doi.org/10.1016/j.jfa.2017.04.008},
}

@article{eskenazis2026functionalperimeterdimensionalbrunnminkowski,
      title={Functional perimeter and the dimensional {B}runn--{B}inkowski inequality for log-concave measures}, 
      author={Alexandros Eskenazis and Apostolos Giannopoulos and Natalia Tziotziou},
      year={2026},
      journal={arXiv:2605.02747} 
}

@article{Armen_Uli,
  title={A new infinitesimal form of the {P}r\'ekopa--{L}eindler inequality with multiplicative structure and applications},
  author={Armeniakos, Sotiris and Ulivelli, Jacopo},
  journal={arXiv:2602.10812},
  year={2026}
}

@article{hos_kol_liv,
  title={On the ${L}_p$-{B}runn--{M}inkowski and dimensional {B}runn--{M}inkowski conjectures for log-concave measures},
  author={Hosle, Johannes and Kolesnikov, Alexander V and Livshyts, Galyna V},
  journal={The Journal of Geometric Analysis},
  volume={31},
  number={6},
  pages={5799--5836},
  year={2021},
  publisher={Springer}
}

@article {Mil_centroaffine,
    AUTHOR = {Milman, Emanuel},
     TITLE = {Centro-affine differential geometry and the log-{M}inkowski
              problem},
   JOURNAL = {J. Eur. Math. Soc. (JEMS)},
  FJOURNAL = {Journal of the European Mathematical Society (JEMS)},
    VOLUME = {27},
      YEAR = {2025},
    NUMBER = {2},
     PAGES = {709--772},
      ISSN = {1435-9855,1435-9863},
   MRCLASS = {58C40 (52A23 52A40)},
  MRNUMBER = {4859577},
MRREVIEWER = {Jacopo\ Ulivelli},
       DOI = {10.4171/jems/1386},
       URL = {https://doi.org/10.4171/jems/1386},
}

@book {Schneider_CB,
    AUTHOR = {Schneider, Rolf},
     TITLE = {Convex bodies: the {B}runn-{M}inkowski theory},
    SERIES = {Encyclopedia of Mathematics and its Applications},
    VOLUME = {151},
   EDITION = {expanded},
 PUBLISHER = {Cambridge University Press, Cambridge},
      YEAR = {2014},
     PAGES = {xxii+736},
      ISBN = {978-1-107-60101-7},
   MRCLASS = {52-02 (52A20 52A39)},
  MRNUMBER = {3155183},
MRREVIEWER = {Andrea\ Colesanti},
}

@article {Ivaki_Milman_p,
    AUTHOR = {Ivaki, Mohammad N. and Milman, Emanuel},
     TITLE = {{$L^p$}-{M}inkowski problem under curvature pinching},
   JOURNAL = {Int. Math. Res. Not. IMRN},
  FJOURNAL = {International Mathematics Research Notices. IMRN},
      YEAR = {2024},
    NUMBER = {10},
     PAGES = {8638--8652},
      ISSN = {1073-7928,1687-0247},
   MRCLASS = {52A40 (53A10)},
  MRNUMBER = {4749179},
MRREVIEWER = {Shi-Zhong\ Du},
       DOI = {10.1093/imrn/rnad319},
       URL = {https://doi.org/10.1093/imrn/rnad319},
}

@article{Ivaki_Hu_weighted_centroaffine,
  title={Weighted centro-affine {P}oincar\'e inequalities},
  author={Hu, Yingxiang and Ivaki, Mohammad N.},
  journal={arXiv:2606.04774},
  year={2026}
}

@article {Cor_Fra_Mau,
    AUTHOR = {Cordero-Erausquin, D. and Fradelizi, M. and Maurey, B.},
     TITLE = {The ({B}) conjecture for the {G}aussian measure of dilates of
              symmetric convex sets and related problems},
   JOURNAL = {J. Funct. Anal.},
  FJOURNAL = {Journal of Functional Analysis},
    VOLUME = {214},
      YEAR = {2004},
    NUMBER = {2},
     PAGES = {410--427},
      ISSN = {0022-1236,1096-0783},
   MRCLASS = {60G15 (60B11 60E15)},
  MRNUMBER = {2083308},
MRREVIEWER = {Ilya\ S.\ Molchanov},
       DOI = {10.1016/j.jfa.2003.12.001},
       URL = {https://doi.org/10.1016/j.jfa.2003.12.001},
}

@article {Esk_Mos_Gauss,
    AUTHOR = {Eskenazis, Alexandros and Moschidis, Georgios},
     TITLE = {The dimensional {B}runn-{M}inkowski inequality in {G}auss
              space},
   JOURNAL = {J. Funct. Anal.},
  FJOURNAL = {Journal of Functional Analysis},
    VOLUME = {280},
      YEAR = {2021},
    NUMBER = {6},
     PAGES = {Paper No. 108914, 19},
      ISSN = {0022-1236,1096-0783},
   MRCLASS = {52A40 (52A20 60E15)},
  MRNUMBER = {4193767},
MRREVIEWER = {Grigory\ M.\ Ivanov},
       DOI = {10.1016/j.jfa.2020.108914},
       URL = {https://doi.org/10.1016/j.jfa.2020.108914},
}

@article {Kol_Liv_GZ,
    AUTHOR = {Kolesnikov, Alexander V. and Livshyts, Galyna V.},
     TITLE = {On the {G}ardner-{Z}vavitch conjecture: symmetry in
              inequalities of {B}runn-{M}inkowski type},
   JOURNAL = {Adv. Math.},
  FJOURNAL = {Advances in Mathematics},
    VOLUME = {384},
      YEAR = {2021},
     PAGES = {Paper No. 107689, 23},
      ISSN = {0001-8708,1090-2082},
   MRCLASS = {52A40},
  MRNUMBER = {4238914},
MRREVIEWER = {Alexandros\ Eskenazis},
       DOI = {10.1016/j.aim.2021.107689},
       URL = {https://doi.org/10.1016/j.aim.2021.107689},
}

@article{boroczkyLYZ2013,
  title={The logarithmic {M}inkowski problem},
  author={B{\"o}r{\"o}czky, K{\'a}roly and Lutwak, Erwin and Yang, Deane and Zhang, Gaoyong},
  journal={Journal of the American Mathematical Society},
  volume={26},
  number={3},
  pages={831--852},
  year={2013}
}

@article{HanZhuPoincare,
  title={{S}harp {P}oincar\'e {I}nterpolation {A}long {W}asserstein {G}eodesics},
  author={Han, Bang-Xian and Zhu, Zhuo-Nan},
  journal={arXiv:2607.10769},
  year={2026}
}

@article{mil_nak_tsu,
  title={{G}aussian inequalities with different covariance matrices- the {H}eterovariant {G}aussian {B}runn-{M}inkowski inequality},
  author={Milman, Emanuel and Nakamura, Shohei and Tsuji, Hiroshi},
  journal={arXiv:2610.13596},
  year={2026}
}

@inproceedings{Latala2002ICM,
  author    = {Rafa{\l} Lata{\l}a},
  title     = {On some inequalities for {G}aussian measures},
  booktitle = {Proceedings of the International Congress of Mathematicians},
  volume    = {II},
  pages     = {813--822},
  year      = {2002},
  publisher = {Higher Education Press},
  address   = {Beijing}
}

@article{yang2026uniformbounddimensionalbrunnminkowski,
      title={A uniform bound in the dimensional {B}runn--{M}inkowski inequality for even log-concave measures}, 
      author={Kai-Wen Yang},
      year={2026},
      journal={arXiv:2607.10104}
}

@article {LocalGlobalPDE,
    AUTHOR = {Chen, Shibing and Huang, Yong and Li, Qi-Rui and Liu, Jiakun},
     TITLE = {The {$L_ p$}-{B}runn-{M}inkowski inequality for {$p<1$}},
   JOURNAL = {Adv. Math.},
  FJOURNAL = {Advances in Mathematics},
    VOLUME = {368},
      YEAR = {2020},
     PAGES = {107166, 21},
      ISSN = {0001-8708,1090-2082},
   MRCLASS = {52A40},
  MRNUMBER = {4088419},
MRREVIEWER = {Julio\ Bernu\'es},
       DOI = {10.1016/j.aim.2020.107166},
       URL = {https://doi.org/10.1016/j.aim.2020.107166},
}

@article {Putter2021,
    AUTHOR = {Putterman, Eli},
     TITLE = {Equivalence of the local and global versions of the
              {$L^p$}-{B}runn-{M}inkowski inequality},
   JOURNAL = {J. Funct. Anal.},
  FJOURNAL = {Journal of Functional Analysis},
    VOLUME = {280},
      YEAR = {2021},
    NUMBER = {9},
     PAGES = {Paper No. 108956, 20},
      ISSN = {0022-1236,1096-0783},
   MRCLASS = {52A39 (52A40)},
  MRNUMBER = {4220744},
MRREVIEWER = {Weidong\ Wang},
       DOI = {10.1016/j.jfa.2021.108956},
       URL = {https://doi.org/10.1016/j.jfa.2021.108956},
}

@article{cor-er_liv_rot,
  title={Remarks on the {B}rascamp--{L}ieb inequality},
  author={Cordero-Erausquin, Dario and Livshyts, Galyna V. and Rotem, Liran},
  journal={arXiv:2610.34289},
  year={2026}
  }

@article {EskeNayTkoczMixtures,
    AUTHOR = {Eskenazis, Alexandros and Nayar, Piotr and Tkocz, Tomasz},
     TITLE = {Gaussian mixtures: entropy and geometric inequalities},
   JOURNAL = {Ann. Probab.},
  FJOURNAL = {The Annals of Probability},
    VOLUME = {46},
      YEAR = {2018},
    NUMBER = {5},
     PAGES = {2908--2945},
      ISSN = {0091-1798,2168-894X},
   MRCLASS = {60E15 (52A20 52A40 60G15 94A17)},
  MRNUMBER = {3846841},
MRREVIEWER = {J.\ H\"usler},
       DOI = {10.1214/17-AOP1242},
       URL = {https://doi.org/10.1214/17-AOP1242},
}

@article{Kol_liv_Rot,
      title={On p-{B}runn-{M}inkowski and {B}rascamp-{L}ieb inequalities}, 
      author={Alexander V. Kolesnikov and Galyna Livshyts and Liran Rotem},
      year={2026},
      journal={arXiv:2507.12099} 
}
\bibliographystyle{siam}

\bigskip

\parbox[t]{9cm}{
Sotiris Armeniakos\\
Institut f\"ur Stochastik und Wirtschaftsmathematik, TU~Wien
e-mail: sotirios.armeniakos@tuwien.ac.at}

\medskip

\parbox[t]{8.5cm}{
Alexander Kolesnikov\\
HSE University, Russian Federation\\
e-mail: sascha77@mail.ru}

\medskip 

\parbox[t]{8.5cm}{
Jacopo Ulivelli\\
Institut f\"ur Diskrete Mathematik und Geometrie, TU Wien
e-mail: jacopo.ulivelli@tuwien.ac.at}
\end{document}